\documentclass{amsart}
\usepackage{amssymb}

\newcommand{\reg}{^{\text{reg}}}
\newcommand{\semi}{^{\text{ss}}}
\newcommand{\mult}{^\times}

\newcommand{\lsup}[1]{{}^{#1}}

\newcommand{\R}{{\mathbb{R}}}
\newcommand{\C}{{\mathbb{C}}}
\newcommand{\Z}{{\mathbb{Z}}}
\newcommand{\N}{{\mathbb{N}}}
\newcommand{\tR}{\widetilde\R}

\newcommand{\bG}{{\mathbf{G}}}
\newcommand{\bH}{{\mathbf{H}}}

\newcommand{\bM}{{\mathbf{M}}}
\newcommand{\bN}{{\mathbf{N}}}
\newcommand{\bP}{{\mathbf{P}}}
\newcommand{\bT}{{\mathbf{T}}}
\newcommand{\bU}{{\mathbf{U}}}
\newcommand{\bS}{{\mathbf{S}}}
\newcommand{\bZ}{{\mathbf{Z}}}
\newcommand{\bA}{{\mathbf{A}}}
\newcommand{\bX}{{\mathbf{X}}}

\let\gg\relax
\newcommand{\gg}{{\mathfrak{g}}}
\newcommand{\bgg}{{\boldsymbol\gg}}
\newcommand{\mm}{{\mathfrak{m}}}
\newcommand{\pp}{{\mathfrak{p}}}
\newcommand{\bb}{{\mathfrak{b}}}
\newcommand{\uu}{{\mathfrak{u}}}

\newcommand{\bmm}{{\boldsymbol\mm}}
\newcommand{\ttt}{{\mathfrak{t}}}
\newcommand{\bttt}{{\boldsymbol\ttt}}
\DeclareMathOperator{\mexp}{\mathsf{e}}

\newcommand{\OO}{{\mathcal{O}}}
\newcommand{\BB}{{\mathcal{B}}}

\newcommand{\checkhat}[1]{\check{\hat{#1}}}
\DeclareMathOperator{\tr}{tr}

\DeclareMathOperator{\res}{res}
\DeclareMathOperator{\ad}{ad}
\DeclareMathOperator{\Ad}{Ad}
\DeclareMathOperator{\Int}{Int}
\DeclareMathOperator{\Gal}{Gal}
\DeclareMathOperator{\Hom}{Hom}
\DeclareMathOperator{\supp}{supp}

\newcommand{\set}[2]{
    {\left\{\left.
    #1\vphantom{#2\bigl(\bigr)}\,\right|
    \,#2\right\}}}

\newcommand{\sett}[2]{\set{#1}{\text{#2}}}

\title  [Local character expansion]
    {The local character expansion near a tame, semisimple element}

\author{Jeffrey D. Adler}
\thanks{The first-named author was partially supported by
the National Security Agency (H98230-05-1-0251)
and by a University of Akron Faculty Research Grant (\#1604).}
\email{adler@uakron.edu}
\address{The University of Akron\\ Akron, OH  44325-4002}

\author{Jonathan Korman}
\email{jkorman@math.toronto.edu}
\address{The University of Toronto\\ Toronto, Ontario M5S 3G3}


\begin{document}

\numberwithin{equation}{section}

\theoremstyle{plain}
\newtheorem{thm}[equation]{Theorem}
\newtheorem{lem}[equation]{Lemma}
\newtheorem{cor}[equation]{Corollary}
\newtheorem{propn}[equation]{Proposition}

\theoremstyle{definition}
\newtheorem{rmk}[equation]{Remark}
\newtheorem{hyp}[equation]{Hypothesis}
\newtheorem{defn}[equation]{Definition}

\begin{abstract}
Consider the character of an irreducible admissible
representation of a $p$-adic reductive group.
The Harish-Chandra-Howe local expansion
expresses this character
near a semisimple element as a linear combination
of Fourier transforms of nilpotent orbital integrals.
Under mild hypotheses, we describe an explicit
region on which the local character expansion is valid.
We assume neither that the group is connected,
nor that the underlying field has characteristic zero.
\end{abstract}

\maketitle

\setcounter{section}{-1}

\section{Introduction}
Let $G$ denote the group of $k$-points of a reductive $k$-group
$\bG$, where $k$ is a nonarchimedean local field.
To simplify the
present discussion, assume for now that $\bG$ is connected and that $k$
has characteristic zero.
Let $(\pi,V)$ denote an irreducible
admissible representation of $G$. Let $dg$ denote a fixed Haar
measure on $G$, and let $C_c^\infty(G)$ denote the space
of complex-valued, locally constant, compactly supported functions on $G$.
The \emph{distribution character} $\Theta_\pi$ of
$\pi$ is the map $C_c^\infty(G) \rightarrow \C$ given by
$\Theta_\pi(f) := \tr \pi(f)$, where $\pi(f)$ is the (finite-rank)
operator on $V$ given by $\pi(f) v := \int_G f(g) \pi(g)v \, dg$.
From Howe~\cite{Howe} and
Harish-Chandra~\cite{Harish-Chandra:submersion}, the distribution
$\Theta_\pi$ is represented by a locally constant function on the
set of regular semisimple elements in $G$.
We will denote this
function also by $\Theta_\pi$.

For any semisimple $\gamma \in G$,
the local character expansion about $\gamma$
(see~\cite{Howe-Fourier} and \cite{Queens})
is the identity
$$
\Theta_\pi (\gamma  \mexp(Y)) = \sum_{\OO} c_{\OO} \hat\mu_\OO
(Y),
$$
valid for all regular semisimple $Y$ in the Lie algebra $\mm$ of
the centralizer of $\gamma$ such that $Y$ is close enough to $0$.
Here, the sum is over the set of nilpotent orbits $\OO$ in $\mm$;
$\hat\mu_\OO$ is the function that represents the distribution
that is the Fourier transform of the orbital integral $\mu_\OO$
associated to $\OO$; $c_{\OO} = c_{\OO,\gamma}(\pi)\in\C$; and
$\mexp$ is the exponential map, or some suitable substitute.

This is a qualitative result, in the sense that it gives no
indication of how close $Y$ must be to $0$ in order for the
identity to be valid.
Many questions in harmonic analysis on $G$ require more quantitative
versions of such qualitative results.
As an example of a quantitative result,
DeBacker~\cite{Homogeneity}
has determined (under some hypotheses on $G$)
a neighborhood of validity for the local
character expansion near the identity,
thus verifying a conjecture of Hales, Moy, and Prasad
(see~\cite{Moy-Prasad-1}).

In this paper, we generalize DeBacker's result for any semisimple
$\gamma\in G$ satisfying mild tameness hypotheses.
(See \S\ref{sec:hyp} for the hypotheses, and
Corollary~\ref{cor:homog-theta} for a precise statement of the
main result.)
When $\gamma$ is regular, we recover a
generalization of Theorem~19 of \cite{loc-const}.

We have taken care not to assume that $k$ has characteristic
zero (though some of our hypotheses do indeed restrict the
characteristic of $k$).
The general philosophy is that
whatever is true in characteristic zero should be true in
large enough positive characteristic.
DeBacker's results
on homogeneity~\cite{Homogeneity} do not assume characteristic zero.
Rather, he keeps careful track of which specific hypotheses
are necessary in order to make each step of his proof valid.
(Previously, the validity of the local character expansion was unknown
in positive characteristic.)
We have done likewise.

We have also taken care not to assume that $\bG$ is connected.
Nonconnected groups are of interest for several reasons,
one of which is that their characters are equivalent to
twisted characters of connected groups.
In order to handle this case, we need to know
that characters are represented by locally
constant functions on the regular set.
All proofs of this result that are known to us assume
either that $k$ has characteristic zero or that $\bG$ is connected.
But it is not difficult to show that the result is true generally.
We include the details in \S\ref{sec:loc-const}.

Note that Theorem~2.1.5(3) of~\cite{Homogeneity}
plays a key role in our proof.

\textbf{Acknowledgements.} We thank Stephen DeBacker, Fiona
Murnaghan, Ju-Lee Kim, and Loren Spice for helpful conversations.
We also thank the referee for helpful comments.

\section{Notation and conventions}
Let $k$ denote a nonarchimedean local field,
and let $\nu$ denote a discrete valuation on $k$. For any
algebraic extension field $E$ of $k$, $\nu$ extends uniquely to a
valuation (also denoted $\nu$) of $E$. Fix a complex-valued,
additive character $\Lambda$ on $k$ that is nontrivial on the ring
$R$ of integers in $k$ and trivial on the prime ideal of $R$.

For a reductive $k$-group $\bG$,
let $\bG^\circ$ denote its connected part,
and
let $\bgg$ denote its Lie algebra.
Let $\bgg^*$ denote the dual of $\bgg$.
Let $G=\bG(k)$, the group of $k$-rational points of $\bG$;
and let $\gg = \bgg(k)$ and $\gg^* = \bgg(k)^*$.
Let $\bZ_\bG$ denote the center of $\bG$.

We use similar notation and font conventions for other groups.
That is, given a group $\bM$, we have
$\bmm$, $M$, etc.

Let $\Ad$ (resp.~$\ad$) denote the adjoint or coadjoint
representation of $\bG$ (resp.~$\bgg$)
on $\bgg$ or $\bgg^*$.
Let $\Int$ denote the conjugation action of $G$ on itself.
For an element or subset $S$ in $G$ and an element or subset
$L$ in $\gg$ or $\gg^*$ (resp.~$G$), we will sometimes write
$\lsup{S} L$ instead of $\Ad(S)L$
(resp.\ $\Int(S)L$).

An element $g\in G$ is \emph{semisimple}
if $\Ad(g)$ is a semisimple linear transformation of $\gg$.
When $g\in G^\circ$, this is equivalent to $g$ belonging
to a torus.
For a subset $S$ of $G$, let $S\semi$ denote the set of semisimple
elements in $S$ (so $S\semi=S\cap G\semi$).
An element
$g\in G$ is \emph{regular semisimple} if the coefficient of $t^n$
in $\det(t-1+\Ad(g))$ is nonzero
(where $n$ is the rank of the component $g\bG^\circ$ in $\bG$;
see~\cite{Clozel}).
We denote the set of regular semisimple elements in $G$ by
$G\reg$.
Similarly we say that an element $X\in\gg$ is
\emph{regular semisimple} if the coefficient of $t^n$ in
$\det(t-\ad(X))$ is nonzero.
We denote the set of regular
semisimple elements in $\gg$ by $\gg\reg$.

For a subset $S$ of $\gg$ (resp.~$G$) let $[S]$ denote
the characteristic function of $S$ on $\gg$ (resp.~$G$).

Call an element $X\in \gg$ \emph{nilpotent}
if there is some one-parameter subgroup $\lambda$ of $\bG$
defined over $k$ such that
$\lim_{t\rightarrow 0} \lsup{\lambda(t)} X = 0$.
Let $\mathcal{N}_\gg$ denote the set of nilpotent elements in $\gg$,
and $\OO_\gg(0)$
the set of nilpotent orbits
under the adjoint action of $G$ on
$\mathcal{N}_\gg$.
We will leave out the subscript when it is understood.
One can similarly define a set $\mathcal{N}^*$ of nilpotent
elements in $\gg^*$.

For any compact group $K$, let $K^\wedge$ denote the set
of equivalence classes of
irreducible, continuous representations of $K$.
We will not always distinguish between a representation
and its equivalence class.
Recall that if $K$ is abelian, then $K^\wedge$ is a group.

Let $\tR := \R \cup \set{ r{+} }{r \in \R}$
and extend the ordering on $\R$ to one on $\tR$ as follows:
for all $r,s\in \R$,
\begin{align*}
r < s{+} & \quad\text{if and only if}\quad  r\leq s; \\
r{+} < s{+}&\quad\text{if and only if}\quad  r< s; \\
r{+} < s{\phantom{+}}& \quad\text{if and only if}\quad  r< s.
\end{align*}
If $r\in\R$, define $(r{+}){+}$ to be $r{+}$.
There is a natural way to extend the additive structure
on $\R$ to an additive structure on $\tR\cup\{\infty\}$.

For any $k$-group $\bG$,
let $\bX^*(\bG)$ denote the lattice of characters of $\bG$.
For any $k$-torus $\bT$ in $\bG$,
let $\Phi=\Phi (\bG,\bT)\subset \bX^*(\bT)$ denote the absolute
root system of $\bG$ with respect to $\bT$.
We can also interpret $\Phi$ as the set of nontrivial eigencharacters
for the adjoint action of $\bT$ on $\bgg$.

Suppose $\bT$ is maximal in $\bG$.
If $\alpha\in\Phi$,
then let $\bU_\alpha$ denote the root group corresponding
to $\alpha$.
This group need not be
defined over $k$, but is defined over a splitting
field $E$ of $\bT$ over $k$.
Let $\bgg(E)_{\alpha}\subset \bgg(E)$ denote the root space
corresponding to $\alpha$.

Let $\Psi(\bG,\bT)$ denote the set of affine roots
of $\bG$ with respect to $\bT$ and $\nu$.
If $\psi \in \Psi(\bG,\bT)$, let
$\dot{\psi}\in \Phi(\bG,\bT)$ denote the gradient of $\psi$.
We denote the root lattice in
$\bgg(E)_{\dot{\psi}}$ corresponding to
$\psi$ by $\bgg(E)_{\psi}$~\cite[3.2]{Moy-Prasad-1}.

\section{Apartments and buildings.}
\label{sec:bldgs}
For any extension $E/k$
of finite ramification degree, let $\BB(\bG,E)$
denote the extended Bruhat-Tits building of $\bG$ over $E$.
Note that if $E/k$ is Galois,
then $\BB(\bG,k)$ embeds naturally in the set
of $\Gal(E/k)$-fixed points of $\BB(\bG,E)$, with equality
when $E/k$ is tame
(see~\cite[(5.11)]{Prasad}).

Every maximal $k$-split torus $\bS$ in $\bG$ has an
associated apartment $\mathcal{A}(\bS,k)$ in
$\BB(\bG,k)$.
Let $\bT$ be a maximal
$k$-torus in $\bG$ containing $\bS$.
Then
$\bT$ splits over some Galois extension $E$, so
$\bT$ has an apartment $\mathcal{A}(\bT,E)$ in
$\BB(\bG,E)$.
The Galois fixed point set of the
apartment of $\bT$ in $\BB(\bG,E)$ is the
apartment of $\bS$ in $\BB(\bG,k)$~\cite[\S2.6]{Tits}.

Suppose $\bM \subset \bG$ is an $E$-Levi $k$-subgroup (that is,
$\bM(E)$ is a Levi subgroup of $\bG(E)$) for some finite Galois
extension $E/k$.
There is a natural family of
$\Gal(E/k)$-equivariant embeddings of $\BB(\bM,E)$ into
$\BB(\bG,E)$.
When $E/k$ is tame, this in turn induces a
natural family of embeddings of $\BB(\bM,k)$ into
$\BB(\bG,k)$.
(In general, there is no canonical way to
pick a distinguished member of this family.
However, all such
embeddings share the same image, and no statement we make will
depend on the choice of embedding.)

We shall require a generalization of the previous paragraph.

\begin{lem}
\label{lem:embedding}
Let $\gamma\in G\semi$ and $\bM=C_{\bG}(\gamma)$.
Then for every tame extension $E/k$ of finite
ramification degree, we have an
$\bM(E)$-equivariant,
toral embedding
${i_E}\colon {\BB(\bM,E)}\longrightarrow {\BB(\bG,E)}$.
Moreover, these maps can be chosen to be compatible in the following sense:
if $E$ and $E'$ are two such fields and $E\supseteq E'$,
then the restriction of $i_E$ to $\BB(\bM,E')$ is $i_{E'}$.
\end{lem}

\begin{proof}
For each Galois extension $E /k$ of finite ramification degree,
Theorem~2.2.1 and Lemma 2.3.4 of~\cite{landvogt} together provide us with a
family $I_{E }$ of toral,
$\bM^\circ(E )$-equivariant embeddings
of $\BB(\bM^\circ,E )=\BB(\bM,E )$ into $\BB(\bG^\circ,E )=\BB(\bG,E )$.
Moreover, $\Gal(E /k)$ acts on $I_{E }$.
For an embedding $f\in I_{E }$ and an element $m\in \bM(E )$,
define $\lsup{m}f\colon \BB(\bM,E ) \rightarrow \BB(\bG,E )$
by
$\lsup{m}f(x) = m f (m^{-1} x)$.
Since $\lsup{m}f = f$ for all $m\in \bM^\circ(E )$, we thus have an
action of the compact group $\bM(E )/\bM^\circ(E )\rtimes\Gal(E /k)$.
By the convexity of $I_{E }$, this action must have
fixed points, which then form a family of toral embeddings
that are equivariant under both $\Gal(E /k)$ and $\bM(E )$.

We now prove the final statement of the lemma.
Let $E/k$ be the maximal tame
subextension of a finite extension over which
both $\bM^\circ$ and $\bG^\circ$ split.
Choose an embedding $i_E$ as in the previous paragraph.
For any tame $L/E$, $\BB(\bM,E)$ generates $\BB(\bM,L)$
as an $\bM(L)$-space, so our choice of $i_E$ determines a choice
of $i_L$.
Thus, it will be enough
to show compatibility between $i_L$ and $i_{L'}$
for all tame, Galois $L/k$ and all $L\supseteq L' \supseteq k$.
The restriction of $i_L$ to $\BB(\bM,L')$ has image
in the set of $\Gal(L/L')$-fixed points of $\BB(\bG,L)$.
But from~\cite{Rousseau} or~\cite{Prasad},
this set is $\BB(\bG,L')$.
\end{proof}

\section{Moy-Prasad filtrations.}
\label{sec:filtrations}
Let $\bT$ be a maximal $k$-torus in $\bG$ with splitting
field $E$ over $k$.
Let $\bT(E)_0$ denote the parahoric subgroup of $\bT(E)$.
For $r\in \tR$, define
$$
\bttt(E)_r :=
\sett{X\in\bttt(E)}{$\nu(d\chi (X))\geq r$ for all $\chi\in\bX^*(\bT)$}
$$
and for $r > 0$,
$$
\bT(E)_r :=
\sett{t\in\bT(E)_0}{$\nu(\chi (t)-1)\geq r$ for all $\chi\in\bX^*(\bT)$}.
$$

For each $(x,r)\in \BB(\bG,E) \times \tR$,
Moy and Prasad
define lattices $\bgg(E)_{x,r}$ in $\bgg(E)$
and $\bgg(E)^*_{x,r}$ in $\bgg(E)^*$.
When $r \geq 0$, they define a normal subgroup $\bG(E)_{x,r}$
of the parahoric subgroup $\bG(E)_x$ of $\bG(E)$.
In particular,
for all $x\in\mathcal{A}(\bT,E)$ and $r\in \tR$,
\begin{equation}
\label{eqn:gxr-defn}
\bgg(E)_{x,r}
=\bttt(E)_r \oplus
\underset{\psi\in\Psi(\bG,\bT),\psi(x)\geq r}\sum \bgg(E)_{\psi}
\end{equation}
Similarly, $\bG(E)_{x,r}$ is defined in terms of the filtrations
on $\bT(E)$ and on root groups.
These definitions depend on the normalization of the valuation $\nu$;
our normalization agrees with that of Yu~\cite{Yu}.
Thus, for example, for any $\alpha\in k\mult$,
$\alpha\cdot \gg_{x,r} = \gg_{x,r+\nu(\alpha)}$.

However, the definitions do not depend on the choice of $\bT$ containing
$x$ in its apartment.
Note that for all $r\in\tR$,
$\bgg(E)_{x,r+} = \cup_{s > r} \bgg(E)_{x,s}$
and $\bG(E)_{x,|r|+} = \cup_{s > |r|} \bG(E)_{x,s}$.

Moy and Prasad also define $\gg_{x,r}$ and $G_{x,r}$
(irrespective of whether or not $\bG$ is $k$-split).
The above normalization
was chosen to have the following property~\cite[1.4.1]{Adler}:
when $E/k$ is tame and $x \in \BB(\bG,k)$,
we have
$$
\gg_{x,r} = \bgg(E)_{x,r} \cap \gg, \qquad
\text{and (for  $r > 0$),}\quad
{G}_{x,r} = \bG(E)_{x,r} \cap G.
$$

We will also use the following notation.
For $r\in \tR$, let
$$
\gg_r=\cup_{x\in\BB(\bG,k)}\gg_{x,r}
\quad\text{and (for $r \geq 0$)}
\quad
G_r=\cup_{x\in\BB(\bG,k)}G_{x,r}.
$$
It is proven in~\cite{Adler-DeBacker} that
$\gg_r$ (resp.~$G_{r}$) is a $G$-domain: a
$G$-invariant, open and closed subset of $\gg$ (resp.~$G$).

For any $x\in\BB(\bG,k)$
and any $0 < r \leq t \leq 2r$, the group $(G_{x,r}/G_{x,t})$
is abelian.
Under many conditions
(for example, if $\bG$ contains a tamely
ramified maximal torus, or if $\bG$ is simply connected),
there exists a ($G_x$-equivariant)
isomorphism (see~\cite{Moy-Prasad-1} or~\cite{Yu})
\begin{equation}
\label{MP-map}
G_{x,r}/G_{x,t} \longleftrightarrow
\gg_{x,r}/\gg_{x,t}\;,
\end{equation}
and thus an isomorphism
\begin{equation}
\label{Kirillov}
(G_{x,r}/G_{x,t})^{\wedge}  \longleftrightarrow
\gg^*_{x,(-t)+}/\gg^*_{x,(-r)+}\;.
\end{equation}
Yu~\cite{Yu:models}
has defined a more complicated filtration on $\bT$ than the one above.
Using this filtration to define $G_{x,r}$, he shows
that \eqref{Kirillov} is valid for all $\bG$.
However, for the groups that we will consider, Yu's
filtration is equivalent to the one above.

\section{Singular depth}
\label{sec:singular}
From now on, $\bG$ is a
reductive $k$-group,
$\gamma \in G\semi$,
and $\bM$ is the centralizer
of $\gamma$ in $\bG$.

\begin{defn}
For $m\in M$, let
$
s(m) = s^{\bG}_{\bM}(m) := \max_{\alpha\in A_m} \{0,\, \nu(\alpha - 1)\},
$
where $A_m$ is the set of generalized eigenvalues of the action of $\Ad(m)$
on $\gg/\mm$.
\end{defn}

\begin{rmk}
If $\bG$ is connected and $m$ is regular,
then the definition of $s(m)$ given above agrees
with the definition
in~\cite[\S1]{loc-const}.
\end{rmk}

\begin{rmk}
\label{regrmk}
Note that $s(m z)=s(m)$ for all $z\in Z_G$
and that $s(h m h^{-1}) = s(m)$ for all $h\in M$.
\end{rmk}

\begin{rmk}
\label{rmk:eigenspaces}
Suppose $m\in M$, and $E/k$ is an extension that contains
all of the generalized eigenvalues of both $\Ad(\gamma)$ and $\Ad(m)$ acting
on $\gg/\mm$.
Since $\Ad(\gamma)$ and $\Ad(m)$ commute,
we can write $\bgg(E)/\bmm(E)$ as a direct sum of subspaces
$V_i$, where $V_i$ is simultaneously an $\alpha_i$-eigenspace
for $\Ad(\gamma)$ and a generalized $\beta_i$-eigenspace for $\Ad(m)$.
\end{rmk}

\begin{lem}
\label{lem:depth-vs-s}
If $m \in M_r\semi$, then $s(m) \geq r$.
\end{lem}

\begin{proof}
Pick a maximal $k$-torus $\bT$ in $\bM$ with $m\in T$,
and a splitting field $E$ for $\bT$.
Then
$$
m\in M_r\cap T \subset \bM(E)_r \cap \bT(E)  = \bT(E)_r,
$$
where the last equality follows from
Theorem~4.1.5(1) of~\cite{DeBacker:group}.
Thus, $\nu (\chi(m) - 1) \geq r$ for all $\chi\in\bX^*(\bT)$.
In particular, this is true for all
$\chi \in \Phi(\bG,\bT)\subset \bX^*(\bT)$.
\end{proof}

\begin{lem}
\label{deepness}
If $m \in M_{s(\gamma)+}$,
then
$s(\gamma m)=s(\gamma)$.
If $\gamma$ is compact mod $Z_G$, then so is $\gamma m$, and conversely.
If $\gamma$ is semisimple, then so is $\gamma m $.
\end{lem}

\begin{proof}
Let $V_i$, $\alpha_i$, and $\beta_i$ be as in Remark~\ref{rmk:eigenspaces}.
For all $i$,
Lemma~\ref{lem:depth-vs-s} implies that
$\nu(\beta_i-1) > s(\gamma)\geq \nu(\alpha_i-1)$, and thus
$$
\nu(\alpha_i\beta_i - 1)
= \nu\bigl((\alpha_i - 1)(\beta_i-1)-(\alpha_i-1)+(\beta_i-1)\bigr)
= \nu(\alpha_i-1).
$$
Thus, $s(\gamma m) = s(\gamma)$.
The second statement of the lemma follows from the fact that
$\alpha_i$ is a unit if and if so is $\alpha_i\beta_i$.
If $\Ad(m)$ is diagonalizable over some field extension,
then $\Ad(m)$ and $\Ad(\gamma)$ are simultaneously diagonalizable,
so the last statement follows.
\end{proof}

For $m \in M$, following~\cite[\S18]{Queens}, define
$$
D_{G/M}(m):=\det\Bigl((\Ad(m)-1)\bigl.\bigr|_{\gg/\mm}\Bigr).
$$
(When $M^\circ =G^\circ$, define $D_{G/M}\equiv 1$.)

For $r\geq 0$ let
\begin{align*}
M'_r  &:= \set{m \in M_r }{ D_{G/M}(\gamma m)\neq 0}, \\
M''_r &:= \set{m \in M_r }{ \gamma m\in G\reg}.
\end{align*}
Note that $M''_r \subset M'_r$,
and these are open, dense subsets of $M_r$.

\begin{cor}
\label{deepness cor}
$M'_{s(\gamma)+}=M_{s(\gamma)+}$.
\end{cor}

\begin{proof}
Let $m\in M_{s(\gamma)+}$.
Let $V_i$, $\alpha_i$, and $\beta_i$ be as in Remark~\ref{rmk:eigenspaces}.
As in the proof of Lemma~\ref{deepness},
$\nu(\alpha_i\beta_i-1) = \nu(\alpha_i-1)$ for all $i$.
Thus,
$$
|D_{G/M}(\gamma m)|
= \prod_i |\alpha_i\beta_i-1|^{\dim V_i}
= \prod_i |\alpha_i-1|^{\dim V_i}
= |D_{G/M}(\gamma)| \neq 0.
$$
Therefore, $m\in M'_{s(\gamma)+}$.
\end{proof}

\section{Intertwining}
\begin{defn}
Let $K$ be a compact open subgroup of $G$ and let $d\in K^{\wedge}$.
For $g\in G$, recall that
$\lsup{g} d$
is the representation of $g K g ^{-1}$ given as
$\lsup{g} d (g k g^{-1}):=d(k)$.
\end{defn}

\begin{defn}
If $d$ and $d'$ are continuous representations of compact
subgroups $K$ and $K'$ (respectively) of $G$, then let $[d:d'] =
\dim_{\C}\Hom_{K\cap K'}(d,d')$.
\end{defn}

\begin{lem}
\label{[d:chi]}
Let $K$ and $L$ be compact subgroups of $G$, and let
$N$ be a compact subgroup of $K$.
Let $d\in K^\wedge$ and let $\chi$ denote a one-dimensional
representation of $L$.
Let $d_1\oplus\cdots\oplus d_n$ be a decomposition of $d$
into representations
of $N$.
If $0\neq [\chi :d]$ then for some $d_i$, $0\neq [\chi :d_i]$.
\end{lem}

\begin{proof}
We have
$0 < [\chi:d] \leq [\chi:d|_N] = \sum_i [\chi:d_i]$.
Therefore, $[\chi:d_i] > 0$ for some $i$.
\end{proof}

\section{Partial traces}
From now on, let $(\pi,V)$ denote
an irreducible admissible representation of $G$.
Let $\Theta_\pi$ denote the distribution character of $\pi$.
This distribution is represented by a
locally constant function (also denoted $\Theta_\pi$) on $G\reg$
(see \S\ref{sec:loc-const}).
Let $\rho(\pi)$
denote the depth of $\pi$~\cite[\S5]{Moy-Prasad-1}.

For any irreducible representation $d$ of a compact open subgroup
$K$, let $V_d$ denote the $(K,d)$-isotypic subspace of $V$.
Let $E_d$ denote the $K$-equivariant projection from $V$ to $V_d$.
Define the distribution $\Theta_d$ by
$\Theta_d(F) := \tr (E_d \,\pi (F)\, E_d)$ for all $F\in C_c^\infty (G)$.
Then $\Theta_d$, which can be thought of as the
`partial trace of $\pi$ with respect to $d$', is represented by
the locally constant function $\Theta_d (x):= \tr (E_d\, \pi (x) \,E_d)$
on $G$.
It follows from the definitions that
\begin{equation}
\label{K-types}
\Theta_\pi (F)=\sum_{d\in K^\wedge}\Theta_d (F)
\quad\text{for all $F\in C_c^\infty (G)$}.
\end{equation}
Note that for each fixed $F$,
$\pi(F)$ has finite rank, so all but finitely many terms in this
sum vanish.

\begin{lem}
\label{conjlem}
\begin{enumerate}
\item
If $x\in G$ and $k\in K$, then
$\Theta_d(kxk^{-1})=\Theta_d(x)$.
\item
If $h \in\ker(d)$,
then $\Theta_d(xh)=\Theta_d(x)$.
\end{enumerate}
\end{lem}

\begin{proof}
The first statement is~\cite[Lemma 14]{loc-const}.
Since $\pi(h)E_d = E_d$, we have
$$
\Theta_d(xh)=\tr (E_d\, \pi (x) \pi (h) \,E_d)
=\tr (E_d \,\pi (x) \,E_d)
=\Theta_d(x).
\qed
$$
\renewcommand{\qed}{}
\end{proof}

Let $N$ be a compact open subgroup of $K$,
and let $d\in K^\wedge$.
Considered as a
representation of $N$, $d$ decomposes into a finite sum of
distinct irreducible representations $d_i$ with multiplicities
$\alpha_i$:
$$
\alpha_1 d_1\oplus\cdots\oplus \alpha_n d_n\,.
$$
For each $i$, let $V_{d,d_i}$ denote the $(N,d_i)$-isotypic
subspace of $V_d$.
Let $E_{d,d_i}$ denote the $N$-equivariant projection
from $V$ to $V_{d,d_i}$.
For $x\in G$,
define
$$
\Theta_{d,d_i}(x) := \tr(E_{d,d_i} \pi(x) E_{d,d_i}).
$$

\begin{rmk}
\label{rmk:theta-ddi}
Note that $\Theta_d(x) = \sum_i \Theta_{d,d_i}(x)$.
Moreover, $\Theta_{d,d_i}$ has invariance properties analogous
to those given for $\Theta_d$ in
Lemma~\ref{conjlem}.
\end{rmk}

\begin{propn}
\label{intertwine}
Let $g\in G$.
If $\Theta_d (g) \neq 0$,
then $0\neq [d : \lsup{g}d]$.
If $\Theta_{d,d_i} (g) \neq 0$,
then $0\neq[d_i : \lsup{g}d_i]$.
\end{propn}

\begin{proof}
Define a pairing $\langle\phantom{x},\phantom{x}\rangle$
on $V_d$ with respect to which $d$ is unitary.
Let $\{v_j\}$ be a basis for $V_d$.
For $k_1,k_2\in K$,
\begin{align*}
\Theta_d(k_1xk_2)
&= \tr(E_d \, \pi(k_1xk_2) \, E_d ) \\
&= \sum_j \langle E_d \pi(k_1 x k_2) v_j, v_j \rangle \\
&= \sum_j \langle \pi(k_1)  E_d \pi(x) \pi(k_2) v_j, v_j \rangle.
\end{align*}
Fixing $k_2$ and letting $k_1$ vary,
we have a sum of matrix coefficients of $d$.
Note that by Lemma~\ref{conjlem}(1)
we have,
\begin{align*}
\Theta_d(k_1xk_2) &= \Theta_d(k_2 k_1 x k_2 k_2^{-1})\\
&=\Theta_d(k_2 k_1 x)\\
&= \sum_j \langle \pi(k_2)  E_d \pi(k_1 x) v_j, v_j \rangle.
\end{align*}
Fixing $k_1$ and letting $k_2$ vary, we again have a sum of matrix
coefficients of $d$.
Therefore, our first statement follows
from Corollary~14.3 of \cite{Queens}.

To prove the second statement, note that
since $E_{d,d_i}$ is a projection onto a subspace of $V_{d_i}$,
if $\Theta_{d,d_i}(g) \neq 0$,
then $\Theta_{d_i}(g) \neq 0$.
\end{proof}

\begin{lem}
\label{lem:int-partial-trace}
Fix $i$.
Let $\chi$ be a character of a closed subgroup $N'$ of $N$.
Let $g\in G$.
Suppose
$$
0 \neq \int_{N'} \Theta_{d,d_i}(gn) \chi(n) \, dn \, .
$$
Then $d_i|_{N'} = \bar\chi$.
\end{lem}

\begin{proof}
For all $n\in N'$, $\pi(n)E_{d,d_i}$ acts on
$V_{d,d_i}$ via the scalar $d_i(n)$.
Let $\{a_j\}$ be the multi-set of diagonal entries
of a matrix that represents, with respect to some basis,
the action of $E_{d,d_i}\pi(g)$ on $V_{d,d_i}$.
Then for some $j$,
$$
0 \neq a_j \int_{N'} d_i(n) \chi(n) \, dn.
$$
The conclusion then follows from Lemma~14.2 of~\cite{Queens}.
\end{proof}

\section{Harmonic analysis}
\subsection*{From distributions on $G$ to distributions on $M$.}
For the distribution $\Theta_\pi$ (or any other distribution) on
$G$, we follow Harish-Chandra~\cite[\S18]{Queens} (see
also~\cite{Rodier}) to define a distribution $\theta$ on $M$ that
captures the behavior of $\Theta_\pi$ near $\gamma$.

Fix $r \geq 0$. The proof of Proposition~1 of~\cite{Rodier} is
valid for general reductive groups,
so we may apply it to see that the following (surjective) map is
everywhere submersive:
\[
\begin{array}{ccccccc}
G \times M_r'
   & \stackrel{p}{\longrightarrow}
   & \lsup{G}(\gamma M_r') \\
(g,m) & \longmapsto & g\gamma m g^{-1}.
\end{array}
\]

\begin{thm}
\label{submersion} There exists a unique, surjective, linear map
\[
\begin{array}{ccccccc}
C_c^\infty (G \times M_r')
   & \longrightarrow
   & C_c^\infty(\lsup{G}(\gamma M_r')) \\
\alpha & \longmapsto & f_\alpha \; ,
\end{array}
\]
such that for all $F\in C_c^\infty (G)$
$$
\int_G F(x)f_\alpha(x) \,dx = \int_{G \times M'_r} F(g\gamma m
g^{-1})\alpha (g,m) \,dg \,dm,
$$
and $\supp(f_\alpha) \subset p (\supp(\alpha))$.
\end{thm}
\begin{proof}
This is Theorem 11, p.~49, of~\cite{Harish-Chandra} applied to the
map $p$ above. (Note that in \emph{loc.\ cit.}, it is not assumed
that the characteristic of $k$ is zero nor that $\bG$ is
connected.)
\end{proof}

\begin{rmk}
The set $\lsup{G}(\gamma M'_r)$ is an \emph {open} (since $p$ is
submersive), $G$-invariant neighborhood of $\gamma$ in $G$, so
$C_c^\infty(\lsup{G}(\gamma M'_r)) \subset C_c^\infty(G)$. Note
that the set $\lsup{G}(\gamma M_r)$ is not necessarily open.
\end{rmk}

Fix an open compact subgroup $K$ of $G$; let $1_K$ denote its
characteristic function. We have the following diagram:
\[
\begin{array}{ccccccccc}
C_c^\infty(M_r)
   & \stackrel{\text{restr.}}{\longrightarrow}
   & C_c^\infty(M_r')
   & \longrightarrow
   & C_c^\infty(G \times M_r')
   & \longrightarrow
   & C_c^\infty(\lsup{G}(\gamma M_r'))
\\
\phantom{.} \\
f & \longmapsto & f & \longmapsto & \alpha & \longmapsto &
f_\alpha
\end{array}
\]
where the first arrow is the restriction map; the second arrow is
given by $\alpha (g,m):=1_K (g)f(m)$; and the third arrow is the
map of Theorem~\ref{submersion}. Note that the support of
$f_\alpha$ is contained in $\lsup{K}(\gamma M'_r)$, an open,
$K$-invariant neighborhood of $\gamma$ in $G$.

Let $\Theta=\Theta_\pi$ denote the distribution character of
$(\pi,V)$.

\begin{defn}
\label{defn:theta}
Define the distribution $\theta$ on $M_r$ by
$\theta (f):=\Theta_\pi (f_\alpha)$.
\end{defn}

\begin{lem}
\label{sum_K-types} Normalize the measure on $G$ so that $K$ has
total measure $1$. Then for each $f\in C_c^\infty (M_r)$
$$
\theta (f)= \sum_{d\in K^\wedge}\int_{M_r}\Theta_d (\gamma m) f(m)
\, dm,
$$
where the sum is over a \emph{finite} set (which depends on $f$)
of representations.
\end{lem}

Note that a similar statement appears on p.~78 of \cite{Queens}.

\begin{proof}
Combine equation \eqref{K-types} with
\begin{align*}
\Theta_d(f_\alpha)
&:=\int_G \Theta_d(x)f_\alpha(x) \,dx \\
&\phantom{:}=
\int_{G \times M'_r}\Theta_d(g\gamma m g^{-1})1_K(g)f(m) \,dg \,dm\\
&\phantom{:}=
    \int_{K \times M'_r} \Theta_d(k \gamma m k^{-1})f(m) \,dk \,dm\\
&\phantom{:}= \int_{M'_r} \Theta_d(\gamma m)f(m) \,dm
\qquad \text{(by Lemma~\ref{conjlem})} \\
&\phantom{:}= \int_{M_r} \Theta_d(\gamma m)f(m) \,dm. \qed
\end{align*}
\renewcommand{\qed}{}
\end{proof}

For the following lemma, note that $\lsup{K}(\gamma M''_r)$ is a
$K$-invariant neighborhood of $\gamma$ in $G\reg$, and that
$\lsup{K}(\gamma M''_r)$ is dense in $\lsup{K}(\gamma M'_r)
\supset \supp(f_\alpha)$. From \S\ref{sec:loc-const}, $\Theta_\pi$ is
represented by a locally constant function on $G\reg$.

\begin{lem}
\label{M''_r} The distribution $\theta$ is represented on $M''_r$
by the function $\theta(m):=\Theta_\pi (\gamma m)$.
\end{lem}
\begin{proof}
For all $f\in C_c^\infty (M_r)$,
the support of
the function $f_{\alpha}$ is a subset of $G\reg$,
so the first integral below makes sense:
\begin{align*}
\theta (f) = \Theta_\pi (f_\alpha)
&:=\int_{G\reg} \Theta_\pi(x)f_\alpha(x) \,dx \\
&\phantom{:}= \int_{\lsup{K} (\gamma M''_r)} \Theta_\pi(x)f_\alpha(x) \,dx\\
&\phantom{:}=
\int_{K \times M''_r} \Theta_\pi (k \gamma m k^{-1})f(m) \,dk \,dm \\
&\phantom{:}= \int_{M''_r} \Theta_\pi (\gamma m) f(m) \,dm. \qed
\end{align*}
\renewcommand{\qed}{}
\end{proof}

\section{Hypotheses}
\label{sec:hyp}
From now on,
we will make certain assumptions:
Let $\bG$, $\gamma$, and $\bM$ be as in \S\ref{sec:singular}.
That is, $\bG$ is a reductive linear algebraic $k$-group,
$\gamma\in G\semi$, and $\bM$ is the centralizer of $\gamma$ in $\bG$.

\begin{hyp}
\label{hyp:tame}
The eigenvalues of $\Ad(\gamma)$ belong to a
tamely ramified extension of $k$.
\end{hyp}

\begin{hyp}
\label{hyp:decomp}
The order of
$\gamma\bZ_\bG\bG^\circ \in \bG/\bZ_\bG\bG^\circ$
is prime to the residue characteristic of $k$.
\end{hyp}

We will use this hypothesis to prove Lemma~\ref{lem:gxr-eigen-decomp}.
However, there are other conditions that would also suffice.

\begin{hyp}
\label{hyp:good-depth}
Using some system of embeddings as in Lemma~\ref{lem:embedding}
to identify
$\BB(\bM,E)$ with a subset of $\BB(\bG,E)$ for each tame $E/k$,
we have that for each such $E$, and
for $r\in \tR$ and $x\in \BB(\bM,k)$,
\begin{align*}
\bmm(E)^*_r &= \bmm(E)^* \cap \bgg(E)^*_r
	& \bmm(E)^*_{x,r} &= \bmm(E)^* \cap \bgg(E)^*_{x,r}  \\
\bM(E)_r &= \bM(E) \cap \bG(E)_r
	& \bM(E)_{x,r}&=\bM(E)\cap \bG(E)_{x,r}  && (r > 0 ).
\end{align*}
\end{hyp}

We will not pursue here the question of when this hypothesis holds.
For now, we note that
it holds when one of the following holds for some tame $E/k$:
$\bM^\circ$ is an $E$-Levi subgroup of $\bG^\circ$;
$\bM^\circ$ and $\bG^\circ$ are both $E$-split, and have the same $E$-split
rank;
or $\bG  =R_{E/k}\bH \rtimes \Gal(E/k)$
for some $E$-group $\bH$
(where $R_{E/k}$ denotes restriction of scalars),
and $\gamma\in \Gal(E/k)$.
(See Lemmata~2.2.3 and~2.2.9 of~\cite{Adler-DeBacker-2}.)

\begin{hyp}
\label{hyp:HB}
There is a nondegenerate
$G$-invariant symmetric bilinear form $B$ on $\gg$ such that we
can identify $\gg^*_{x,r}$ with $\gg_{x,r}$ via the map
$\gg\rightarrow\gg^*$ defined by $X \longmapsto (Y \mapsto B(X,Y))$.
\end{hyp}

(Groups satisfying this hypothesis
are discussed in~\cite[\S4]{Adler-Roche}.)
Thus, we can (and eventually will)
identify $\gg_{x,r}$ with $\gg^*_{x,r}$
and $\mm_{x,r}$ with $\mm^*_{x,r}$

The following hypothesis concerns
the existence of a ``mock exponential'' map.
\begin{hyp}
\label{hyp:HM}
Let $e = \max\{\rho(\pi),s(\gamma)\}$.
There exists a
bijection $\mexp:\gg_{e+}\longrightarrow G_{e+}$ such that
for all $x \in \BB(\bG,k)$ and
all $ e<r \leq s\leq t \leq 2s$,
we have that
$\mexp$ induces the group isomorphism
$\gg_{x,s}/\gg_{x,t} \rightarrow G_{x,s}/G_{x,t}$
of \eqref{MP-map}.
Moreover,
for all $s>e$,
$\mexp(\mm_s) = M_s$, and
\begin{enumerate}
\item
\label{item:cosets}
for $r,s>e$,
all $X\in \mm_{x,r}$,
and all $Y\in \mm_{x,s}$, we have
$\mexp(X)\cdot\mexp(Y)\equiv \mexp(X+Y)$
modulo $M_{x,(r+s)}$;
\item
for all $m\in M$ we have
$\Int(m)\circ(\mexp|_{\mm_{e+}}) =(\mexp|_{\mm_{e+}}) \circ \Ad(m)$;
\item
\label{item:adjoint}
for all $x\in \BB(\bM,k)$, all $s,t\in\tR$ with $s>e$,
all $Y\in \gg_{x,s}$,
and all $X \in \gg_{x,t}$,
we have
$\lsup{\mexp(Y)} X - X \in \gg_{x,s+t}$.
\end{enumerate}
\end{hyp}

Note that item (\ref{item:adjoint}) in the hypothesis asserts,
for all $x\in\BB(\bM,k)$,
a weaker version of Proposition~1.6.3
of~\cite{Adler}, and the remainder of the hypothesis
is a weaker version of Hypothesis~3.2.1 in
\cite{Homogeneity}.
Item (\ref{item:cosets}) implies that
$\mexp$ carries a Haar measure on $\mm$ into a Haar measure on $M$.

Mock exponential maps are known to exist in several situations.
For example, for $\mathbf{GL}_n(k)$, the map $X \mapsto (1+X)$ works.
For a classical group that splits over a tame extension of
$k$, with odd residue characteristic, the Cayley transform works.
If $k$ has characteristic zero and $e$ is sufficiently large,
then the exponential map works.

We will need the next two hypotheses
in order to apply Theorem~\ref{homog_thm}.

\begin{hyp}
\label{hyp:char0}
Assume that $\bM$ satisfies Hypothesis~3.4.3 of \cite{Homogeneity},
concerning the convergence of nilpotent orbital integrals.
\end{hyp}
This is automatically satisfied when $k$ has characteristic zero.

\begin{hyp}
\label{hyp:homog}
Assume that $k$ and $\bM$ satisfy
the hypotheses of Theorem~2.1.5(3) of \cite{Homogeneity}.
\end{hyp}

\section{Lie algebra decompositions}
\label{sec:decomp}

\begin{lem}
\label{lem:gxr-eigen-decomp}
Let $E/k$ be a tame extension containing the eigenvalues of $\Ad(\gamma)$.
For all $x\in \BB(\bM,k)$, and $r\in \tR$,
we have
$$
\bgg(E)_{x,r} = \bigoplus_\alpha(\bgg(E)_{x,r}\cap \bgg(E)_{\gamma,\alpha}) \, ,
$$
where the sum is over the set of eigenvalues of the action
of $\Ad(\gamma)$ on $\bgg(E)$, and each $\bgg(E)_{\gamma,\alpha}$
is the associated eigenspace.
\end{lem}

\begin{proof}
Let $n$ be the order of
$\gamma\bZ_\bG\bG^\circ \in \bG/\bZ_\bG\bG^\circ$.
Then
${\gamma}^n \in (\bZ_\bG)(E) \bT(E)$ for some
maximal torus $\bT$ such that $\bT/\bZ_\bG^\circ$ is $E$-split.
Since ${\gamma}^n$ acts trivially on $\bttt(E)$ and acts on each $\bT$-root
space in $\bgg(E)$,
we have from \eqref{eqn:gxr-defn} that
$$
\bgg(E)_{x,r} = \bigoplus_\beta(\bgg(E)_{x,r}\cap \bgg(E)_{{\gamma}^n,\beta}),
$$
where the sum runs over the eigenvalues $\beta$ of $\Ad({\gamma}^n)$,
and each $\bgg(E)_{{\gamma}^n,\beta}$ denotes the corresponding eigenspace.
Thus it is enough to show that for each $\beta$,
$$
\bgg(E)_{x,r}\cap \bgg(E)_{{\gamma}^n,\beta}
=\bigoplus_{\alpha^n = \beta} (\bgg(E)_{x,r} \cap \bgg(E)_{\gamma,\alpha}).
$$
This will follow from Lemma~\ref{lem:lattice-decomp}
and Hypothesis~\ref{hyp:decomp}.
\end{proof}

We may identify $\bmm$ with the
$1$-eigenspace of $\Ad(\gamma)$ in $\bgg$.
Define $\bmm^\bot$ to
be the sum of all of the other eigenspaces of $\Ad(\gamma)$.
Identify $\bmm^*$ with $\set{X\in \bgg^*}{X(\bmm^\bot) = 0}$ and
define ${\bmm^*}^\bot := \set{X \in \bgg^*}{X(\bmm) = 0}$.
These objects are all defined over $k$, and
\begin{equation*}
\gg=\mm\oplus\mm^\bot
\quad\text{and}\quad
\gg^*=\mm^*\oplus{\mm^*}^\bot
\end{equation*}

The following result is well known when $\bG$ is connected.
See~\cite[Proposition~1.9.2]{Adler}.

\begin{cor}
\label{cor:gxr-decomp}
Suppose $x \in \BB(\bM,k)$ and $r\in \tR$.
Then
$$
\gg_{x,r} = (\mm \cap \gg_{x,r}) \oplus (\mm^\bot\cap \gg_{x,r}).
\qed
$$
\end{cor}

\begin{lem}
\label{lem:lattice-decomp}
Let $R$ denote
the ring of integers in $k$,
and let $\varpi\in R$ be a uniformizer.
Let $V$ be a finite-dimensional vector space over $k$,
and $L$ a lattice in $V$.
Let $T\colon V\rightarrow V$
denote a diagonalizable linear map such that $T(L) = \varpi^rL$
for some $r\in\Z$.
Let $\alpha_1,\ldots,\alpha_\ell$ be the eigenvalues
of $T$, and $V_1,\ldots,V_\ell$
the corresponding eigenspaces.
Suppose that for $i\neq j$,
we have that $\alpha_i\not\equiv\alpha_j\bmod \varpi^{r+1}R$.
Then
$$
L = \bigoplus_{i=1}^\ell (L \cap V_i).
$$
\end{lem}

\begin{proof}
Without loss of generality, we may replace $T$ by $\varpi^{-r} T$,
and thus assume that $T(L) = L$,
so that $\alpha_i\in R\mult$ for all $i$.
Let $x\in L$, and write $x = \sum x_i$,
where $x_i \in V_i$.
We must show that $x_i\in L$ for all $i$.
Let $S = \set{i}{x_i \notin L}$,
and suppose that $S$ is nonempty.
Pick a minimal nonempty subset $I\subset S$
such that $\sum_{i\in I}\beta_i x_i \in L$
for some coefficients $\beta_i \in R\mult$.
Pick such coefficients, and let $y$ be the resulting sum.
We must have $|I|>1$.
Let $i_0\in I$.
Then $\alpha_{i_0}^{-1} T(y) - y \in L$.
But we may rewrite this element as
$$
\sum_{i\in I \smallsetminus \{i_0\}}
\bigl(\frac{\alpha_i}{\alpha_{i_0}} - 1\bigr) x_i,
$$
contradicting the minimality of $I$.
\end{proof}

\section{Some lemmata}
\label{sec:lemmas}
Some results in this section are stated in terms of the Lie
algebras $\mm$ and $\gg$.
However, by Hypothesis~\ref{hyp:HB},
the analogous
results for $\mm^*$ and $\gg^*$ are also valid (with the same proofs).

\begin{lem}
\label{Z-gamma}
Let $t\in \R$, $x\in\BB(\bM,k)$,
and $\gamma'\in \gamma M_{x,s(\gamma)+}$.
If $Z\in \mm^\bot \cap (\gg_{x,-t} \smallsetminus \gg_{x,(-t)+})$
then $\lsup{\gamma'} Z -Z \not\in \gg_{x,(-t+s(\gamma'))+}$.
\end{lem}

\begin{proof}
Write $\gamma' = \gamma m$, with $m\in M_{x,s(\gamma)+}$.
Let $Z' = \lsup{m} Z - Z$.

From Hypothesis~\ref{hyp:tame}, there is a tame extension $E$ of
$k$ containing the eigenvalues of $\Ad(\gamma)$.
Write
$Z=\sum_{\alpha} Z_\alpha$,
and
$Z'=\sum_{\alpha} Z'_\alpha$,
where each sum is over the set of
eigenvalues for the action of $\Ad(\gamma)$ on $\bgg(E)$, and each
$Z_\alpha$ and $Z'_\alpha$
belongs to the corresponding eigenspace.
Then
$\lsup{\gamma'}Z - Z = \sum_\alpha Y_\alpha$,
where
$Y_\alpha = \alpha Z'_\alpha + (\alpha-1)Z_\alpha$.
From our hypothesis on $Z$, there is some $\alpha$ so that
$Z_\alpha\notin \bgg(E)_{x,(-t)+}$.
From Lemma~\ref{deepness}, $s(\gamma')=s(\gamma)$.
Therefore, from Lemma~\ref{lem:gxr-eigen-decomp}
it will be enough to show that
$Y_\alpha \notin \bgg(E)_{x,(-t+s(\gamma))+}$.
Note that from
Hypothesis~\ref{hyp:HM}(\ref{item:adjoint}),
$Z' \in \gg_{x,(-t+s(\gamma))+}$.
Therefore, Lemma~\ref{lem:gxr-eigen-decomp}
implies that
$Z'_\alpha\in \bgg(E)_{x,(-t+s(\gamma))+}$.

Suppose $\nu(\alpha)\geq 0$.
Then $\alpha Z'_\alpha \in \bgg(E)_{x,(-t+s(\gamma))+}$.
By the definition of $s(\gamma)$,
$(\alpha - 1)Z_\alpha\notin \bgg(E)_{x,(-t+s(\gamma))+}$,
and our conclusion follows.

Now suppose that $\nu(\alpha)<0$.
Let $q= \nu(\alpha) = \nu(\alpha-1)$.
Then
$\alpha Z'_\alpha \in \bgg(E)_{x,(q-t+s(\gamma))+}$
and
$(\alpha-1) Z_\alpha \notin \bgg(E)_{x,(q-t)+}
\supset \bgg(E)_{x,(q-t+s(\gamma))+}$,
so
$Y_\alpha \notin \bgg(E)_{x,(q-t+s(\gamma))+}\supset
\bgg(E)_{x,(-t+s(\gamma))+}$.
\end{proof}

\begin{propn}
\label{proposition}
Let $r\in \R$, $x\in\BB(\bM,k)$,
$X\in \gg$,
and $\gamma'\in \gamma M_{x,s(\gamma)+}$.
Write
$X=Y+Z$ according to the decomposition in
Corollary~\ref{cor:gxr-decomp}.
If ${^{\gamma'} X}-X \in \gg_{x,(-r)+}$,
then $Z\in \gg_{x,(-r-s(\gamma))+}$.
\end{propn}

\begin{proof}
For some $t \in \R$,
$Z\in \mm^\perp\cap(\gg_{x,t}\smallsetminus \gg_{x,t+})$.
By Lemma~\ref{Z-gamma},
$\lsup{\gamma'} Z-Z \not\in \gg_{x,(t+s(\gamma'))+}$.
On the other hand, since $\lsup{\gamma'} X-X \in \gg_{x,(-r)+}$
decomposes as
$(\lsup{\gamma'} Y-Y)+(\lsup{\gamma'} Z-Z)
\in \mm_{x,(-r)+}\oplus (\mm^\perp \cap \gg_{x,(-r)+})$,
we have $\lsup{\gamma'} Z-Z \in \gg_{x,(-r)+}$.
Thus $t+s(\gamma') > -r$,
which implies that $Z\in\gg_{x,t}\subseteq \gg_{x,(-r-s(\gamma'))+}$.
Finally, recall from Lemma~\ref{deepness} that $s(\gamma') = s(\gamma)$.
\end{proof}

\begin{propn}
\label{nilpotent_descent}
If $x\in\BB(\bM,k)$,
$r\in \R$, and $Y\in \mm_{x,r}\smallsetminus \mm_{x,r+}$,
then
$$
(Y+\gg_{x,r+})\cap
\mathcal{N}\neq \emptyset\quad \text{if and only if} \quad
(Y+\mm_{x,r+})\cap \mathcal{N}_\mm\neq
\emptyset.
$$
\end{propn}

\begin{proof}
This follows from
Corollary~3.2.6 of~\cite{Adler-DeBacker}
and Hypothesis~\ref{hyp:good-depth}.
\end{proof}

\begin{defn}
Let $0<r\leq t \leq 2r$.
A character $d\in (G_{x,r}/G_{x,t})^{\wedge}$ is called
\emph{degenerate} if the coset that corresponds to $d$
under the isomorphism \eqref{Kirillov}
contains nilpotent elements.
One can similarly define what it means for a character of
$(M_{x,r}/M_{x,t})$ to be degenerate.
\end{defn}

\begin{lem}
\label{rest_rmk}
Let $x\in\BB(\bM,k)$, $q,t\in\tR$, and suppose $0<q\leq t \leq 2q$.
Let
$d\in (G_{x,q}/G_{x,t})^\wedge$ be a character such that
$[d:\lsup{\gamma'}d]\neq 0$ for some $\gamma ' \in \gamma M_{x,s(\gamma)+}$.
Suppose that for some $v > q+s(\gamma)$,
$d$ is trivial on $M_{x,v}$.
Then $d$ is trivial on $G_{x,v}$.
\end{lem}

\begin{proof}
By \eqref{Kirillov},
$d$ corresponds to some coset
$\Upsilon\in\gg^*_{x,(-t+)} / \gg^*_{x,(-q)+}$.
By Lemma~1.8.1 of~\cite{Adler},
$\Upsilon \cap \lsup{\gamma'}\Upsilon \neq \emptyset$.
Pick $X\in \Upsilon$ such that
$\lsup{\gamma'}X \in \Upsilon \cap \lsup{\gamma'}\Upsilon$.
Then $\lsup{\gamma'}X - X \in \gg^*_{x,(-q)+}$.
Write $X=Y+Z$ with respect
to the decomposition in
Corollary~\ref{cor:gxr-decomp}.
By
Proposition~\ref{proposition},
$Z\in \gg^*_{x,(-q-s(\gamma))+} \subseteq \gg^*_{x,(-v)+}$.
Hence, $X + \gg^*_{x,(-v)+} = Y + \gg^*_{x,(-v)+}$.
Since $d$ is trivial on $M_{x,v}$, we have that $Y\in\mm^*_{x,(-v)+}$.
Therefore $X\in \gg^*_{x,(-v)+}$, implying that $d$ is trivial
on $G_{x,v}$.
\end{proof}

\begin{rmk}
\label{rmk:M-concentrated}
In fact, the proof shows that
$d$ is trivial on a slightly larger subgroup
that corresponds to the lattice
$\mm_{x,s^+} + (\mm^\bot\cap \gg_{x,s})$
(for some $s$ such that $\mm_{x,s^+} = \mm_{x,v}$)
via \eqref{MP-map}.
\end{rmk}

\section{A result of DeBacker}

We recall the following from, for example,
\S4.1 of
\cite{Adler-DeBacker}.
Let
$dX$ be a Haar measure on $\mm$.
For any $C_c^\infty
(\mm)$, we define the \emph{Fourier transform} $\hat{f}\in
C_c^\infty (\mm^*)$ of $f$ by
$$
\hat{f}(\chi)=\int_{\mm} f(X)\cdot \Lambda (\chi (X)) \,dX
$$
for $\chi\in\mm^*$.
Let $d\chi$ be a Haar measure on $\mm^*$.
For $f\in C_c^\infty (\mm^*)$ we use
the natural identification of $\mm^{**}$ with
$\mm$ and define the \emph{Fourier transform} $\hat{f}\in
C_c^\infty (\mm)$ by
$$
\hat{f}(X)=\int_{\mm^*} f(\chi)\cdot \Lambda (\chi (X)) \,d\chi
$$
for $X\in \mm$.
We normalize the measures $dX$
and $d\chi$ so that for $X\in \gg$ and $f\in C_c^\infty (\mm)$
\begin{equation}
\label{eqn:fourier_inversion}
\hat{\hat{f}}(X)=f(-X).
\end{equation}

Recall that, from Hypothesis~\ref{hyp:HB},
we can (and eventually will) identify
$\mm$ with its linear dual $\mm^*$.
With this identification, the Fourier transform becomes a map from
$C_c^\infty (\mm)$ to itself.
Given our normalization of measures, we have that
for $x\in \BB(\bM,k)$ and $r\in \tR$,
$[\mm_{x,r}]^\wedge = [\mm^*_{x,(-r)+}]$.

Using $r$ and $\mexp$ as in Hypothesis~\ref{hyp:HM},
and the distribution $\theta$ from Definition~\ref{defn:theta},
we now define a distribution $\widehat{\theta}$
on $\mm^*$ (compare~\cite[\S2.1]{Homogeneity}).

\begin{defn}
\label{defn:hattheta}
For $f\in C_c^\infty (\mm^*)$, define
$\widehat{\theta}(f):=
\theta (\widehat{f}\,|_{\mm_r} \circ \mexp^{-1})$.
\end{defn}

\begin{rmk}
For $x\in \BB(\bM,k)$,
$f\in C({\mm^*_{x,-s}}/{\mm^*_{x,(-r)+}})$
if and only if
$\widehat{f}\in C({\mm_{x,r}}/{\mm_{x,s+}})\subset C_c^\infty (\mm_r)$.
Hence in this case
$\widehat{\theta}(f)=\theta (\widehat{f} \circ \mexp^{-1})$.
\end{rmk}

The following were defined in~\cite[\S2.1]{Homogeneity}.

\begin{defn}
For any $M$-invariant subset $S$ of $\mm$ or $\mm^*$,
let $J(S)$ denote the space of $M$-invariant distributions
supported on $S$.
\end{defn}

\begin{defn}
\label{J-tilde}
Suppose $-s \leq -r$.
Define the spaces of distributions
$$
\widetilde{J}(\mm)_{x,-s,(-r)+}
=
\sett{T\in J(\mm) }{
    \begin{tabular}{l}
    for  $f\in C(\mm_{x,-s}/\mm_{x,(-r)+})$, if $T(f) \neq 0$ \\
    then $\supp(f)\cap (\mathcal{N}_\mm+\mm_{x,(-s)+}) \neq\emptyset$
    \end{tabular}
    },
$$
and
$$
\widetilde{J}(\mm)_{(-r)+}
=
\underset{x\in\BB(\bM,k)} \bigcap
\underset{-s\leq -r }\bigcap
\widetilde{J}(\mm)_{x,-s,(-r)+}.
$$
\end{defn}

\begin{defn}
Define the space of functions
$$
\mathcal{D}(\mm)_{(-r)+}:=\underset{x\in\BB(\bM,k)}\sum
C_c(\mm/\mm_{x,(-r)+}).
$$
\end{defn}

\begin{rmk}
One can define $\widetilde J(\mm^*)_{x,-s,(-r)^+}$,
$\widetilde J(\mm^*)_{(-r)+}$,
and $\mathcal{D}(\mm^*)_{(-r)+}$ similarly.
\end{rmk}

\begin{rmk}
\label{Fourier_domain}
The function
$f\in C_c^\infty(\mm^*)$ lies in
$C_c^\infty({\mm^*_{r}})$ if and only if
$\widehat{f}\in \mathcal{D}(\mm)_{(-r)+}$
(see Definition 4.2.1 and Lemma 4.2.3 of~\cite{Adler-DeBacker}).
\end{rmk}

If $T$ is a
distribution on $\mm$, then we let
$\res_{\mathcal{D}(\mm)_{(-r)+}}T$ denote the restriction of $T$ to
$\mathcal{D}(\mm)_{(-r)+}$.
The following homogeneity result was proved by
DeBacker~\cite[Theorem 2.1.5]{Homogeneity}.
From a remark in the introduction to \emph{loc.\ cit.},
the result does not require that $\bM$ be connected.

\begin{thm}
\label{homog_thm}
$
\res_{\mathcal{D}(\mm)_{(-r)+}}\widetilde{J}(\mm)_{(-r)+}=
\res_{\mathcal{D}(\mm)_{(-r)+}}
J(\mathcal{N}_\mm) .
$
\end{thm}

\section{The Main Theorem}

Recall that we are assuming the hypotheses in \S\ref{sec:hyp}. Let
$\widehat\theta$ be as in definition~\ref{defn:hattheta}.

\begin{thm}
\label{main_thm}
Let $r > \max\{\rho (\pi), 2s(\gamma)\}$.
Then $\widehat{\theta} \in \widetilde{J}(\mm^*)_{(-r)+}$.
\end{thm}

\begin{proof}
It is enough to show that
$\widehat{\theta} \in \widetilde{J}(\mm^*)_{x,-s,(-r)+}$
for all $x\in \BB(\bM,k)$
and all $s,r\in \R$ such that $s \geq r$.
Fix $x\in \BB(\bM,k)$ and $s \geq r$, and take
$f\in C({\mm^*_{x,-s}}/{\mm^*_{x,(-r)+}})$.
Suppose
$0\neq\widehat{\theta}(f)$.
We need to show that
$\supp(f)\cap (\mathcal{N}^*_\mm+\mm^*_{x,(-s)+}) \neq\emptyset$.
By the linearity
of $\widehat{\theta}$, it suffices to show this for
$f=[Y+\mm^*_{x,(-r)+}]$, where $Y\in \mm^*_{x,-s}$.
In other words, it suffices to show that the character
$\chi := \mexp^{-1} \circ \Lambda_Y$
of $M_{x,s}/M_{x,s+}$
is degenerate,
where $\Lambda_Y:=\Lambda \circ Y$.
We have
$$
\widehat{f}(W)=
\Lambda_Y(W)[\mm_{x,r}](W)=\Lambda_Y(W)\underset{\overline{Z}\in
{\mm_{x,r}}/{\mm_{x,s}} }\sum
[\overline{Z}](W).
$$
Thus $0\neq\widehat{\theta}(f)=\theta (\widehat{f}\circ
\mexp^{-1})$ implies that for some $Z\in \mm_{x,r}$,
$0\neq \theta(h)$,
where
$h\in C_c^\infty (M_r)$ is defined by
$h(\mexp(W))=\Lambda_Y(W)[Z+\mm_{x,s}](W)$.
Using Lemma~\ref{sum_K-types}
with $r':=r-s(\gamma)$,
we get
$$
0\neq \theta (h)
= \sum_{d\in G_{x,r'}^\wedge}\int_{M_r}\Theta_d (\gamma m) h(m) \,dm .
$$
Thus for some $d\in G_{x,r'}^\wedge$,
\begin{equation}\label{int_h}
0\neq
\int_{M_r}\Theta_d (\gamma m) h(m) \, dm
=
\int_{M_{x,r}}\Theta_d (\gamma m) h(m) \, dm,
\end{equation}
where the last equality holds because $h$ is supported on $M_{x,r}$.

Pick minimal $t\in\R$ such that $d$ is trivial on $G_{x,t+}$.
Let $q = \max\{r',(\frac{t}{2})+\}$.
Since $q\geq r'$, one can restrict $d$ to $G_{x,q}$.
Since the
group $G_{x,q}/G_{x,t+}$ is abelian, this restriction
decomposes into a finite sum of irreducible, one-dimensional
representations $d_1\oplus\cdots\oplus d_n$.
From equation
\eqref{int_h} and Remark~\ref{rmk:theta-ddi}, we
see that for some $i$, and for all $v>s$,
\begin{align*}
0& \neq \int_{M_{x,r}}\Theta_{d,d_i} (\gamma m) h(m) \, dm \\
&=\sum_{\dot{n}\in M_{x,r}/M_{x,v}}
    \int_{M_{x,v}} \Theta_{d,d_i}(\gamma \dot{n}m) h(\dot{n}m) \, dm\\
&=\sum_{\dot{n}\in M_{x,r}/M_{x,v}} h(\dot{n})
\int_{M_{x,v}} \Theta_{d,d_i}(\gamma \dot{n}m) \, dm,
\end{align*}
where we think of each sum as running over a set of coset representatives
for $M_{x,r}/M_{x,v}$,
and in the last line above we use the fact that $h$ is
constant on $M_{x,s+}$.
Therefore, for some $\dot{n}$,
\begin{equation}
\label{eqn:nonzero-int-dotn}
0\neq \int_{M_{x,v}}
\Theta_{d,d_i}(\gamma \dot{n}m) \, dm.
\end{equation}
If we assume that $v\geq q$, then we may apply
Lemma~\ref{lem:int-partial-trace}
to see that $d_i$ is trivial on $M_{x,v}$.
Moreover,
\eqref{eqn:nonzero-int-dotn}
implies that
$\Theta_{d,d_i}(\gamma')\neq 0$ for
some $\gamma ' \in \gamma M_{x,s(\gamma)+}$, which by
Proposition~\ref{intertwine} means that
$0\neq [d_i :\lsup{\gamma'} d_i]$.
If we further assume that $v>q+s(\gamma)$, then we may
apply Lemma~\ref{rest_rmk} to see that
$d_i|_{G_{x,v}} \equiv 1$.
Since $d$ is an irreducible
representation of $G_{x,r'}$, it follows that $G_{x,r'}$ permutes
the $d_i$'s transitively.
Therefore, $d$ is trivial on $G_{x,v}$,
so $v > t$.
Since this is true for all $v$ satisfying $v>s$ and
$v> q + s(\gamma)$,
we see that
\begin{equation}
\label{t<=max}
0< t \leq \max\{s,q+s(\gamma)\} \, .
\end{equation}

We now use this inequality to prove four others:
\begin{align}
\label{s>t/2+}
s & > (\textstyle\frac{t}{2})+ \\
\label{r+s>t}
r + s & > t \\
\label{s>=q}
s & >  q \\
\label{s>=t}
s & \geq t\,.
\end{align}
The first two of these follow trivially from the last.
However, we prove them independently because we want to isolate
the one part of this paper, the proof of \eqref{s>=t},
that relies on the hypothesis
that $r>2s(\gamma)$.

Note that \eqref{s>t/2+} is obvious in the case
where $t\leq s$.
So assume that $t\leq q+s(\gamma)$.
Recall that $q$ is either $r'$ or $(\frac{t}{2})+$.
If $q=r'$, then $s\geq r = q+s(\gamma) \geq t > t/2+$.
If $q=(\frac{t}{2})+$, then
$t\leq (\frac{t}{2}{+}) + s(\gamma)$.
Since $t\in\R$,
$t\leq t/2+s(\gamma)$,
so $t/2 \leq s(\gamma) < r \leq s$.
Since $s\in\R$,
we have \eqref{s>t/2+}.

To prove \eqref{r+s>t}, note that \eqref{s>t/2+}
implies that $s(\gamma) + s > (\frac{t}{2}{+}) + s(\gamma)$,
so $r+s> (\frac{t}{2}{+}) + s(\gamma)$.
Since $r+s>r$, we have $r+s>\max\{r,\frac{t}{2}+s(\gamma){+}\}
=q+s(\gamma)$.
Since $r+s>s$, we have $r+s > \max\{s,q+s(\gamma)\} \geq t$.

To prove \eqref{s>=q}, note that $s\geq r > r'$.
From \eqref{s>t/2+}, we conclude that
$s>\max\{r',\frac{t}{2}{+}\} = q$.

To prove \eqref{s>=t}, we use \eqref{t<=max}
to reduce to the case where $q+s(\gamma) \geq t$.
If $q=r'$, then
$s\geq r = r'+s(\gamma) = q + s(\gamma) \geq t$.
If $q = \frac{t}{2}{+}$,
then since $q+s(\gamma) \geq t$,
we have $s(\gamma){+} \geq \frac{t}{2}$.
Since $s(\gamma)\in\R$,
$s(\gamma) \geq \frac{t}{2}{+}$.
Therefore (using the hypothesis that $r>2s(\gamma)$), we have
$$
s \geq r \geq 2s(\gamma){+} = (s(\gamma){+}) + s(\gamma)
\geq (\textstyle\frac{t}{2}{+}) + s(\gamma)
= q+s(\gamma) \geq t.
$$

We have
\begin{align*}
0 &\neq \int_{M_{x,r}}\Theta_{d,d_i} (\gamma m) h(m) \, dm \\
&= \int_{\mm_{x,r}}
        \Theta_{d,d_i} (\gamma \mexp(W)) \Lambda_Y(W)[Z+\mm_{x,s}](W)
        \,dW\\
&=\Lambda_Y(Z)\int_{\mm_{x,s}}
    \Theta_{d,d_i}(\gamma \mexp(Z+W)) \Lambda_Y(W)\,dW.
\end{align*}
Now let $z:=\mexp(Z)$.
From \eqref{r+s>t}, $d_i$ is trivial on $M_{x,r+s}$.
Apply
Hypothesis~\ref{hyp:HM}(\ref{item:cosets}) together with
Remark~\ref{rmk:theta-ddi} to obtain:
\begin{equation*}
0 \neq \int_{\mm_{x,s}}
    \Theta_{d,d_i} (\gamma z \mexp(W)) \Lambda_Y(W) \,dW\\
= \int_{M_{x,s}}
    \Theta_{d,d_i} (\gamma z m) \chi (m) \,dm.
\end{equation*}
From \eqref{s>=q} and Lemma~\ref{lem:int-partial-trace},
$0\neq [\bar\chi:d_i]$.
From \eqref{s>=t}, $d_i$ is trivial on $G_{x,s+}$.
By Remark~\ref{rmk:M-concentrated},
the restriction $\tilde{d}$ of $d_i$ to $G_{x,s}$ is represented
by a coset $Y+ \gg^*_{x,(-s)+}$ where $Y \in \mm^*_{x,-s}$.
Since $s > \rho(\pi)$,
Theorem~3.5 of~\cite{Moy-Prasad-2} implies that
$\tilde{d}$
is degenerate.
Thus,
$(Y + \gg^*_{x,(-s)+})\cap \mathcal{N}^*\neq \emptyset$.
Use Proposition~\ref{nilpotent_descent} to
conclude that
$(Y +\mm^*_{x,(-s)+})\cap \mathcal{N}^*_\mm\neq\emptyset$,
and hence that $\bar\chi = d_i|_{M_{x,s}}=\tilde{d}|_{M_{x,s}}$ is degenerate.
Thus, $\chi$, is degenerate.
\end{proof}

From now on, use Hypothesis~\ref{hyp:HB} to identify $\mm^*_{x,r}$
with $\mm_{x,r}$ for all $x\in\BB(\bM,k)$ and all $r\in \tR$. For
$\OO\in\OO_\mm(0)$, let $\mu_\OO$ denote the corresponding
nilpotent orbital integral and let $\widehat{\mu}_\OO$ denote its
Fourier transform (both are distributions). Note that in general
the set $\OO_\mm(0)$ can have infinite cardinality. However,
Hypothesis~\ref{hyp:homog} puts restrictions on $k$ and
$\bG$, which guarantee that the cardinality is finite.

\begin{cor}
\label{cor:loc-char-exp}
Let $r > \max\{\rho (\pi), 2s(\gamma)\}$.
Then
$$
\theta (f \circ \mexp^{-1})
= \underset{\OO\in \OO_\mm(0)}\sum
c_\OO\widehat{\mu}_\OO(f) \quad\text{for all $f\in C_c^\infty (\mm\reg_r)$},
$$
where
$c_\OO=c_{\OO,\gamma}(\pi)$ are complex
constants that depend on $\OO$, $\gamma$ and $\pi$.
\end{cor}

\begin{proof}
Let $f\in C_c^\infty(M_r)$.
Then by Remark~\ref{Fourier_domain},
$\widehat{f}\in\mathcal{D}(\mm)_{(-r)+}$.
Let $\checkhat{f}(X) = \widehat{f}(-X)$.
For all $\OO\in\OO_\mm(0)$,
let $-\OO$ denote $\set{X\in\mm}{-X\in\OO}$,
and note that $-\OO\in\OO_\mm(0)$.
Then for some coefficients $\set{c_\OO}{\OO\in\OO_\mm(0)}$,
we have
\begin{align*}
\theta (f \circ \mexp^{-1})
&= \widehat{\theta}(\checkhat{f})
    &&\text{by \eqref{eqn:fourier_inversion}}\\
&=\underset{\OO\in\OO_\mm(0)}\sum c_{-\OO}\mu_{-\OO}(\checkhat{f})
    && \text{by Theorems~\ref{homog_thm} and~\ref{main_thm}} \\
&=\underset{\OO\in \OO_\mm(0)}\sum c_\OO\mu_\OO(\widehat{f}) \\
&=\underset{\OO\in \OO_\mm(0)}\sum c_\OO\widehat{\mu}_\OO(f).
&&\qed
\end{align*}
\renewcommand{\qed}{}
\end{proof}

From Hypothesis~\ref{hyp:char0} and work of Huntsinger (see
Theorem A.1.2 of \cite{Adler-DeBacker-2}), it is known that
$\widehat{\mu}_\OO$ is represented by a locally constant function
(which we will also denote $\widehat{\mu}_\OO$) on $\mm\reg$.
(When $k$ has characteristic zero, this is a result of
Harish-Chandra~\cite[Theorem 4.4]{Queens}.)
From \S\ref{sec:loc-const}, $\Theta_\pi$ is also
represented by a locally constant function.

\begin{cor}
\label{cor:homog-theta}
Let $r > \max\{\rho (\pi), 2s(\gamma)\}$.
Then
$$
\Theta_\pi (\gamma\mexp(Y)) = \underset{\OO\in
\OO_\mm(0)}\sum c_\OO\widehat{\mu}_\OO(Y)
$$
for all $Y\in \mm''_r := \mexp^{-1} (M''_r)$.
\end{cor}

\begin{proof}
For $m\in M_r''$,
we have
$\gamma m \in G\reg$.
Let $\bT = C_\bG(\gamma m)^\circ$.
From Lemmata~1 and 2 of~\cite{Clozel},
$0\neq \det ((\Ad(\gamma m) - 1)_{\gg/\ttt})$.
Since $\Ad(\gamma m) = \Ad(\gamma) \Ad(m) = \Ad(m)$
on $\mm/\ttt$,
$0\neq \det((\Ad(m)-1)|_{\mm/\ttt})$,
and so
$m\in M\reg$.
Thus, $M_r'' \subset M\reg$,
and thus $\mm_r'' \subset \mm\reg$, and so the right-hand
side of the equation makes sense.
Since $\gamma M_r'' \subset G\reg$, the left-hand side makes sense
from Proposition~\ref{propn:loc-const}.
The result now follows from
Corollary~\ref{cor:loc-char-exp} and
Lemma~\ref{M''_r}.
\end{proof}

\begin{cor}
Let $r > \max\{\rho (\pi), 2s(\gamma)\}$.
Then
$$
\Theta_\pi (g\gamma\mexp(Y)g^{-1})=
\underset{\OO\in \OO_\mm(0)}\sum c_\OO\widehat{\mu}_\OO(Y)
\quad\text{for all}\quad
Y\in \mm''_r, \;g\in G.
$$
\end{cor}
\begin{proof}
This follows from Corollary~\ref{cor:homog-theta}
and the $G$-invariance of $\Theta_\pi$.
\end{proof}
\begin{rmk}
When $\gamma$ is regular, we have that $\bM^\circ$ is a torus,
and
so the only nilpotent orbit in $\mm$ is the $0$ orbit.
Thus in this case
there is only one orbital integral in the character expansion and
its Fourier transform is a constant function.
This means that
the domain of validity of the local character expansion near a
regular semisimple element is a domain on which
$\Theta_\pi$ is constant.
Thus we recover a generalization of the main result of~\cite{loc-const}.
\end{rmk}

\begin{rmk}
It would be desirable, for applications of motivic integration to character theory,
to have a version of Theorem~\ref{main_thm}
(and thus of its corollaries) that is valid under the weaker hypothesis
that $r>\max\{\rho(\pi),s(\gamma)\}$.
In order to obtain such a theorem, one would have to replace the
one part of the proof of Theorem~\ref{main_thm}
that assumes $r>2s(\gamma)$:
the proof of inequality \eqref{s>=t}.
This inequality allows us to apply
Remark~\ref{rmk:M-concentrated},
a slight strengthening of Lemma~\ref{rest_rmk},
to the character $d_i$.
However, if we had a version of Lemma~\ref{rest_rmk} strong enough
to apply directly to the representation $d$ (which is not necessarily
one dimensional), then \eqref{s>=t} would be unnecessary.
We will pursue this matter elsewhere.
\end{rmk}

\section{Appendix: Characters are locally constant}
\label{sec:loc-const}
In this section only, $k$ denotes a nonarchimedean local field,
$\bG$ denotes an arbitrary reductive $k$-group,
and $\Theta_\pi$ denotes the distribution character of
an admissible irreducible representation $\pi$ of $G$.
We place no other hypotheses on $k$ or $G$.
In particular, we do not assume any of the hypotheses in
\S\ref{sec:hyp}.

\begin{propn}
\label{propn:loc-const}
The character $\Theta_\pi$ is represented by a function that
is locally constant on $G\reg$.
\end{propn}

Recall that we needed this in order to prove Cor.~\ref{cor:homog-theta}.

Of course, when $k$ has characteristic zero
and $\bG$ is connected, the above is a result of Howe~\cite{Howe} and Harish-Chandra~\cite{Queens}.
(Moreover, Harish-Chandra~\cite{Queens}
showed that $|D_G|^{1/2}\Theta_\pi$ is locally integrable,
a result that Clozel~\cite{Clozel} extended to nonconnected groups.
But we do not need this.)
When $k$ has arbitrary characteristic and $\bG$ is connected,
Harish-Chandra~\cite{Harish-Chandra:submersion}
gave an alternative proof of local constancy (which is independent of local integrability), except that he proves Theorem~1 of \emph{loc. cit.}
(the ``submersion principle'')
only in characteristic zero.
According to Harish-Chandra, Borel knew a characteristic-free
proof of the submersion principle.
However, the first such published proof is due to Gopal Prasad,
and appears as Appendix~B in~\cite{Adler-DeBacker-2}.

Thus, the only novelty in our result above is that we
simultaneously relax the requirements that $k$ have characteristic
zero and that $\bG$ be connected.

The proof, on the other hand, involves no novelty at all,
since the characteristic-free proofs
mentioned above
can be generalized without difficulty to the case of
nonconnected groups.  Here we outline how to do so.

We look first at Prasad's proof of the submersion principle.
Describing the changes necessary to generalize the proof
would take about as much space as the proof itself.
So for the reader's convenience, we present the generalized
proof, even though it is almost a word-for-word copy of Prasad's.

\begin{propn}[Harish-Chandra]
Suppose $k$ is an arbitrary field, $\bG$ is a reductive
$k$-group, and $P$ is
a parabolic subgroup of $G^\circ$.
For $x \in G$, let $x^*$ denote the image of $x$ under the
projection of $G$ onto $G^* = G/P$.
Fix $\gamma \in G\reg$.  The mapping
$$x \rightarrow (\lsup{x}\gamma)^*$$
from $G$ to $G^*$ is everywhere submersive.
\end{propn}

\begin{proof}[Proof (G.~Prasad).]
Let $\phi_\gamma \colon G \rightarrow G^*$ denote the map
$x \rightarrow (\lsup{x}\gamma)^*$.
Since $\phi_\gamma(xy) =
\phi_{\lsup{y}\gamma}(x)$, it is enough to show that $\phi_\gamma$ is
submersive at the identity.

Define
$$V := \{ \lsup{\gamma^{-1}}X - X \, | \, X \in \gg \}.$$
We must show that  $\gg = V + \pp$ where $\pp$ is the Lie algebra of $P$.

Without loss of generality, we assume that $k$ is algebraically
closed.  

The identity component of $C_G(\gamma)$ is a
torus $T$ of $G$. 
Let $B$ be a Borel subgroup containing $T$ and let $U$ denote its
unipotent radical.
Let $T'$ be a maximal torus lying in $B \cap P$.
Let $\bb$ (resp., $\mathfrak{t}$, resp., $\mathfrak{t}'$, resp. $\uu$)
denote the Lie algebra of $B$ (resp. $T$, resp., $T'$, resp. $U$). 
Note that $\gg$ decomposes as a direct sum of $\Ad(\gamma)$-eigenspaces,
and $V$ contains all of the non-trivial eigenspaces,
so $\mathfrak{t} + V = \gg$.
Thus, it is enough to show that $\mathfrak{t} \subset V + \pp$.
Since $\mathfrak{t}'$ is contained in $\pp$ and $\uu\subset V$,
we conclude that $\uu + \mathfrak{t}' \subset V + \pp$.
However, $\bb = \uu + \mathfrak{t}'$,
so  $\mathfrak{t} \subset \bb \subset V + \pp$,
and so we are finished.
\end{proof}

We now turn to the rest of Harish-Chandra's proof of
local constancy.  Now that the submersion principle has
been established in the generality we require,
the rest of the proof works without change,
except for the places where Harish-Chandra uses
the Cartan decomposition in the proof of
his Theorem~2 (\S4 of \cite{Harish-Chandra:submersion}).
We now describe how to generalize this decomposition in a way
that meets the requirements of Harish-Chandra's proof.

Let $K$ denote a special maximal compact subgroup of $G^\circ$,
i.e., the stabilizer of a special vertex $x\in\BB(\bG^\circ,k)$.
Pick a maximal $k$-split torus $\bA$ such that $x$
lies in the apartment of $A = \bA(k)$.
Let $\bM = C_\bG(\bA)$.
Then $\bM$ is a minimal Levi $k$-subgroup of $\bG^\circ$,
and there is some minimal parabolic
$k$-subgroup $\bP$ with Levi decomposition $\bP=\bM\bN$.

Recall that $G^\circ$ has the Cartan decomposition for $G^\circ = K M^+ K$,
where
we can define $M^+$ to be the set of all $m\in M$
such that $m^{-1}(P \cap G_{x,r})m \subseteq P\cap G_{x,r}$
for all $r>0$.
Harish-Chandra uses the Cartan decomposition, together with the following fact:
there is a base $\{K_r\}_{r\in\N}$ of compact, open, normal subgroups of $K$
such that for each $r$, there is a compact open subgroup $P_r\subset P$
such that $m^{-1} P_r m \subset K_r$ for all $m\in M^+$.
So we only need to prove the following:

\begin{lem}
There is a subset $\widetilde{M}^+ \subset G$
such that
$G = K \widetilde{M}^+ K$, and such that
for all $m\in \widetilde{M}^+$, we have that
\begin{equation*}
\tag{$*$}
m^{-1} (P \cap G_{x,r})m  \subseteq G_{x,r} \quad\text{for all $r>0$}.
\end{equation*}
\end{lem}

\begin{proof}
Let $\widetilde{M}$ denote the set of elements of $G$ that normalize
the pair $(P,A)$.  Then $\widetilde{M}\cap G^\circ = M$.
Since $G^\circ$ acts transitively on the set of all pairs consisting
of a maximal $k$-split torus and a minimal parabolic subgroup containing
it, we have that every coset in $G/G^\circ$ contains an element of
$\widetilde{M}$.  That is, $\widetilde{M}/M$ is naturally isomorphic to
$G/G^\circ$.  Thus, $G = K \widetilde{M} K$.

Define
$\widetilde{M}^+$ to be the set of all elements of $\widetilde{M}$
that satisfy condition~($*$) of the lemma.

Since $M/A$ is compact, $M$ has a unique Moy-Prasad filtration:
for $r>0$, $M_r = M_{x,r} = M \cap G_{x,r}$.
Since $\widetilde{M}$ normalizes $M$, it also normalizes each $M_r$.
Each group $G_{x,r}$ has an Iwahori factorization with respect
to $P=MN$.
Thus, an element $m$ of $\widetilde{M}$ lies in $\widetilde{M}^+$
if and only if $m^{-1}(N \cap G_{x,r})m \subseteq N\cap G_{x,r}$
for all $r>0$.

Recall that $A \cap M^+$ contains an element $m_0$ with the property that
$m_0^{-1}(N \cap G_{x,r})m_0 \subset N\cap G_{x,r+1}$
for all $r>0$.
Thus,
for any $m\in\widetilde{M}$,
there is some positive $i$ such that
$m m_0^i \in\widetilde{M}^+$.
This implies that
$\widetilde{M}^+$ meets every coset of $\widetilde{M}/M$,
and thus that $G = K \widetilde{M}^+ K$.
\end{proof}

\end{document}